\documentclass[12pt]{article}

\usepackage{setspace}  \usepackage[margin=1in]{geometry}
\usepackage{amsfonts}
\usepackage{graphicx}
\usepackage{amsmath}
\usepackage{xcolor}
\usepackage{lineno}
\usepackage{hyperref}
\usepackage{bm}
\usepackage{ulem}

\usepackage{caption}
\usepackage{subcaption}

\title{Physics Informed Deep Learning for Strain Gradient Continuum Plasticity}
\author{Ankit Tyagi$^{*}$, Uttam Suman$^{*}$, Mariya Mamajiwala$^{\dagger}$, and Debasish Roy$^{*\#}$  \\
	\small $^{*}$Computational Mechanics Lab, Department of Civil Engineering, Indian Institute of Science, Bangalore, India \\
    \small $^{\dagger}$Department of Computer Science, University of Sheffield, United Kingdom \\
	\small $^{\#}$Centre of Excellence in Advanced Mechanics of Materials, Indian Institute of Science, Bangalore, India}

\date{}

\begin{document}
\maketitle

\begin{abstract}
	We use a space-time discretization based on physics informed deep learning (PIDL) to approximate solutions of a class of rate-dependent strain gradient plasticity models. The differential equation governing the plastic flow, the so-called microforce balance for this class of yield-free plasticity models, is very stiff, often leading to numerical corruption and a
	consequent lack of accuracy or convergence by finite
	element (FE) methods. Indeed, setting up the discretized framework, especially
	with an elaborate meshing around the propagating plastic bands whose locations are often unknown a-priori,
	also scales up the computational effort significantly. Taking inspiration from physics informed
	neural networks, we modify the loss function of a PIDL model in several novel
	ways to account for the balance laws, either through energetics or via the
	resulting PDEs once a variational scheme is applied, and the constitutive
	equations. The initial and the boundary conditions may either be imposed strictly by
	encoding them within the PIDL architecture, or enforced weakly as a part of the loss function. The flexibility in the implementation
	of a PIDL technique often makes for its ready interface with powerful
	optimization schemes, and this in turn provides for many possibilities in
	posing the problem. We have used freely available open-source libraries that
	perform fast, parallel computations on GPUs. Using numerical
	illustrations, we demonstrate how PIDL methods could address the computational challenges posed by strain gradient plasticity models. Also, PIDL
	methods offer abundant potentialities, vis-\'a-vis a somewhat straitjacketed 
	and poorer approximant of FE methods, in customizing the formulation as per the problem objective. 

\end{abstract}

\section{Introduction}
Finite element (FE) methods are the de facto numerical schemes to solve an initial boundary-value problem (IBVP) in solid mechanics. They have many advantages; for example, they are known to solve geometrically complex problems, applying  boundary conditions is straightforward, and the easily constructed $C^0$ FE interpolants are good enough for many problems~\cite{Hughes1987}.  There is also a scope for higher order theories, e.g. mixed approaches~\cite{JamunKumar2022, dhas2022mixed}, to address myriad numerical issues such as locking.

On the flip side, FE methods do suffer from many shortcomings. They are in general computationally expensive and do not scale well with modern accelerated hardware such as graphical processing units (GPUs)~\cite{Pfaff2021}. The mesh generation itself could be a bottleneck, especially in the presence of some local ill-conditioning. Constructing higher order theories is always a challenge. Going beyond $C^0$ finite element interpolation is generally not possible and results in complicated schemes~\cite{Hughes1987}. With greater demands on design, many problems in solid mechanics are now multiscale, and existing computational tools such as the FE struggle to solve them satisfactorily. Among these are problems of damage propagation and fatigue life prediction~\cite{Roy2017,Srinivasa2022}. With FE methods, every problem must be solved from scratch, howsoever incrementally different it might be from an already-solved problem, and there is no scope for extrapolation. That is why they are not suitable for sensitivity analysis and inverse problems. There are also issues related to element distortion and various kinds of locking~\cite{RAJENDRAN20101044, reddy1997locking}. The solution scheme is incremental and repetitive, requiring costly matrix inversion at every step. This leads to all sorts of convergence issues. From an implementation point of view, every novel theory requires a new user element and/or material subroutine, which is quite cumbersome.

Regardless of the great commercial successes of the FE, it is time we looked for alternatives to keep pace with the ever increasing demands on component design. A potential alternative is physics informed deep learning (PIDL), which offers smoother approximants and their derivatives that are dense in the space of square integrable functions~\cite{Hornik1989}. This aspect in itself might ameliorate several shortcomings of FE methods. PIDL methods exploit modern deep learning architectures to solve forward and inverse IBVPs~\cite{Raissi2019}. The solution scheme is not incremental or repetitive and, depending on the optimizer used for training, it does not involve any matrix inversion. The field variables are smoothly approximated with neural networks~\cite{Nielsen2015}. These are typically made up of a successive composition of affine transformations and a non-linear function, called the activation function. Neural networks are universal function approximators~\cite{Hornik1989}, and with a proper choice of the activation function, they are $C^\infty$. Their derivatives may be exactly calculated using automatic differentiation~\cite{Baydin2018}. The solution is achieved by requiring the optimizer to drive the loss function to zero, which corresponds to residues of the field equations going to zero. PIDL has a streamlined implementation procedure, which does not typically require special treatment for specific problems. There is a scope to obtain the solution only over a part of the problem domain, rather than being forced to solve for all the variables all the time.  One may work with the strong form~\cite{Raissi2019} or weak form~\cite{Weinan2018}. Optimized neural networks have been shown to extrapolate and generalize over geometry and boundary conditions~\cite{Pfaff2021}. One benefits from the already solved problems to solve new ones via transfer learning~\cite{Yosinski2014}. PIDL methods are naturally suited to sensitivity analysis and inverse problems~\cite{Haghighat2021}. However, one must bear in mind that the working of the PIDL is generally as good 
as that of the optimizer for the chosen loss function.

Recently PIDL methods have been used to solve many forward problems in mechanics~\cite{Raissi2019, Weinan2018, Samaniego2020}, even though most of these could have been solved as well with FE methods. In this work, we use PIDL to solve a particular class of strain gradient plasticity (SGP) models~\cite{Anand2005, Lele2008a, Lele2008}, on which FE methods particularly struggle. The yield-free SGP theories, unlike classical continuum plasticity models that exploit a yield surface~\cite{Simo1998}, use energetic (and sometimes dissipative) length scales typically associated with terms involving the gradient of the accumulated plastic strain in the free energy in line with certain meso-scale level (`microstructural') observations on the plastic flow. Accordingly, unlike the classical models, they can simulate the well-observed size effect, e.g. in metals undergoing inhomogeneous plastic flow at micron-scale ~\cite{Hutchinson2000, fleck1994strain, stolken1998microbend}. Unlike the classical plasticity models, SGP-based simulations are also mesh-size independent. An FE implementation of the SGP model requires creating a special finite element (e.g. in the form of a user element subroutine in case one uses a commercial software) and a material model (in the form of a user material subroutine). However, the SGP model contains an extremely stiff differential equation, called the microforce balance, to evolve the plastic strain. This is owing to the nature of the constitutive (hardening) equations for the micro-stresses and the presence of length scales (which are quite small). In addition, there is the usual linear momentum balance, also called the macroforce balance, to evolve deformation. Macro- and micro-force balances are generally coupled when plastic flow occurs. The FE solution of the SGP model suffers from numerous convergence issues, and one may have to take ad-hoc measures vis-à-vis the time stepping (e.g. avoiding sharp transition regions of the stress-strain curves) and iterative updates (e.g. manipulating the step-size) to get to a reasonable solution. Although all our numerical experiments are done on as SGP model, most of the observations are equally valid for a phase field model~\cite{miehe2010phase}, which has a somewhat similar formulation.

The rest of the article is organized as follows. Sections 2 and 3 give brief introductions to the PIDL and the SGP model respectively. Section 4 describes the implementation details and the numerical results. Finally, concluding remarks are given in section 5. 

\section{Physics informed deep learning (PIDL)}
\begin{figure}
	\centering
	\includegraphics[width=0.6\textwidth]{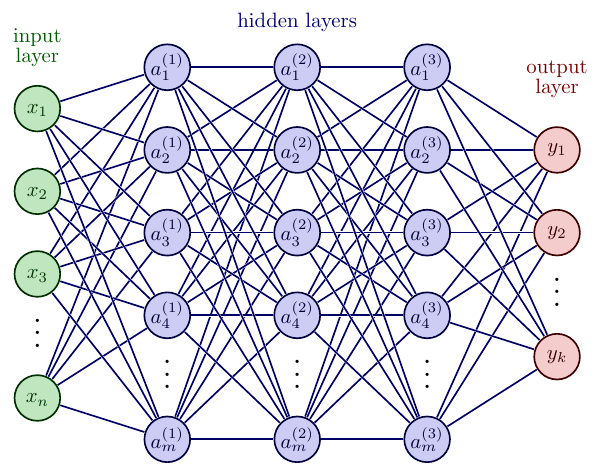}
	\caption{A deep feedforward network (DFN)~\cite{tikz_net}.}
	\label{fig:nn}
\end{figure}

\begin{figure}
	\centering
	\includegraphics[width=0.8\textwidth]{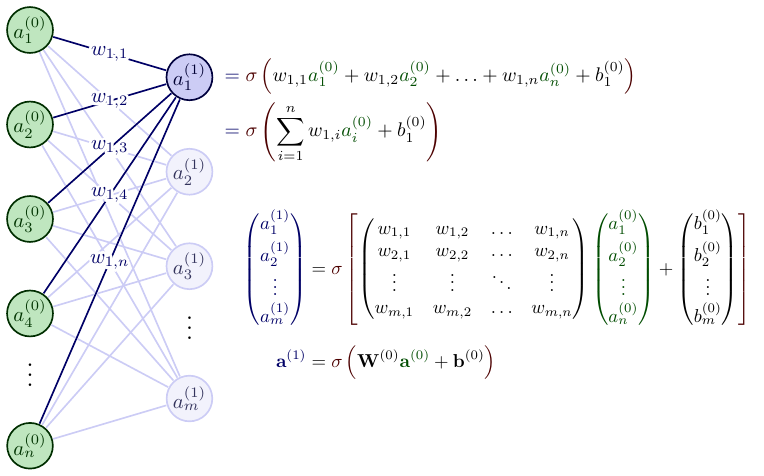}
	\caption{Activation function in one neuron and one layer in matrix notation~\cite{tikz_net}.}
	\label{fig:activation}
\end{figure}

PIDL methods use neural networks to approximate the field variables and deep learning methods to find the parameters (weights and biases) appearing in the approximants. Neural networks are of many kinds, e.g., feedforward, convolutional, recurrent etc. Presently, we use feedforward networks, which are generally employed for regression and classification tasks. A feedforward network is a universal function approximator~\cite{Hornik1989} and has the structure of a parameterized composition function $f:\mathbb{R}^n \rightarrow \mathbb{R}^m$ such that
\begin{gather}
	f(x) = l^{(d)} \circ \sigma \circ l^{(d-1)} \cdots \sigma \circ l^{(1)} \circ \sigma \circ l^{(0)}(x)\ , \\
	\text{where } l^{(i)}(x) = \mathbf{W}^{(i)} x + \mathbf{b}^{(i)}\ .
\end{gather}
$\mathbf{W}$'s and $\mathbf{b}$'s are the neural network parameters, respectively called weights and biases. $\sigma$ is a nonlinear function called the
activation function, and its output is called the activation $a$. $d$ represents the depth of the neural network.
Feedforward networks with $d > 1$ are called deep feedforward networks (DFNs).
Because of its compositional nature, a DFN can be represented graphically as
shown in figure~\ref{fig:nn}, where each node is called a neuron and $d$ denotes the number of hidden layers. For the details of operations that are performed at each neuron and layer, see figure~\ref{fig:activation}.

Traditional deep learning is data-driven, that is, we have a set of input-output pairs and  we want to find the map between them. A neural network is used as an approximant to this map and a loss function is constructed whose minimization determines the network parameters. The loss function is some measure of the gap between the network approximation of the outputs and the actual outputs for the given set of inputs. A common choice of the loss function in regression tasks is the mean squared error (MSE). An optimization algorithm must be used drive the loss value towards zero by adjusting the network parameters. Generally
a first-order gradient-based optimization algorithm, such as the gradient descent with an adjustible learning rate,
is used. The exact gradients of the loss function with respect to the network
parameters are easily obtained using backpropagation, which is a particular implementation of the automatic differentiation~\cite{Baydin2018}. The optimization
process is called the training of the network, and as the network trains, it is
said that the network is learning. Deep learning is the field concerned with
algorithms and techniques that make learning of deep neural networks (DNNs)
possible.

In the field mechanics of materials, we generally do not have access to enough data to perform data-driven deep learning. PIDL goes beyond the realm of data-driven deep learning by using the physics of the problem within the loss function alongside any known data. In PIDL, neural networks are used to approximate the unknown field variables appearing in the energetics or the balance laws. In other words, the governing field equations are incorporated in the loss function in such a way that a
minimization of the latter corresponds to the minimization of the residue of the former. Generally, the loss function is written as a weighted sum of mean squared residues of all field equations for a set of collocation points, also called the training points. The boundary and
initial conditions are either imposed weakly by including them in the loss
function or imposed strictly by encoding them into the architecture of the
neural network. A loss function in PIDL may thus be written as:
\begin{equation}
	\mathcal{L}_{total} = \mathcal{L}_{data} + \mathcal{L}_{physics},
\end{equation}
In this work, in the absence any available experimental data, we take $\mathcal{L}_{data}$ as a measure of the residue in boundary and
initial conditions. $\mathcal{L}_{physics}$ represents the residue in the governing
equations or balance laws. Indeed, the data loss can also incorporate other known data, e.g., some sensor
measurements if available. Overall, PIDL is a very general and modular discretization plus approximation framework, which provides flexibility at every step of problem solving.

\section{Strain gradient plasticity (SGP) theories}
Several well-known experiments, such as torsion of thin wires and indentation tests, have shown that metals at micron-scale undergoing inhomogeneous plastic flow exhibit size effect, which means that the flow strength increases with decreasing size~\cite{Fleck1997}. Conventional continuum plasticity theories cannot simulate such a physical phenomenon, as they do not incorporate any material length-scales. The SGP theories incorporate the material length-scales; thus, they can simulate the size effect, and their simulations should, in principle, be mesh-size independent. The earliest known SGP theory was proposed by Aifantis~\cite{Aifantis1987}, who made the flow stress ($\tau$) not only a function of the accumulated plastic strain ($\gamma$) but also of its gradient. This was done by simply adding a non-local term $-c\nabla^2 \gamma$ to the conventional flow rule:
\begin{equation}
    \tau = \tau(\gamma) - c\nabla^2 \gamma,
\end{equation}
 where $\tau(\gamma)$ denotes the classical flow resistance, $c$ is a positive material constant and $\nabla^2$ represents the Laplace operator. Note that, given $\tau$, the new flow rule has the form of a partial differential equation (PDE) in $\gamma$. It is also called the microscopic force balance~\cite{Gurtin2010}. Subsequently, many such theories have been developed and a class of them are based on the principle of virtual power~\cite{Fleck1997, Fleck2001, Gurtin2004, Gudmundson2004, Gurtin2005}. Among these, we consider the theory given by Gurtin and Anand~\cite{Gurtin2005}, where the relevant gradient field is the plastic strain tensor, and the resulting flow rule is a tensor-valued PDE. This makes a theory far richer than that of Aifantis, where the relevant gradient field is the scalar accumulated plastic strain, and the flow rule is also scalar-valued. The theory is capable of accounting for backstress, an energy dependent on the Burgers vector, and strengthening~\cite{Gurtin2010}. Standard finite element implementations of the theory have been attempted~\cite{Anand2005, Lele2008a, Lele2008}. Typically, plastic strain is added to the nodal degrees of freedom. Thanks to the constitutive relations, a large number of terms appear in the Jacobian of the residue. Many of these terms are however simply neglected to reduce the computation cost and, perhaps more importantly, make the model less stiff. There have also been attempts to use mixed FE methods, but these are either implemented using a yield surface~\cite{miehe2013mixed} or for only strain gradient elasticity theories~\cite{choi2023mixed}.

In this work, we attempt at solving 1D and 2D models of a strip of infinite length and finite width ($0\leq y \leq h$) undergoing plastic deformation under simple shear~\cite{Anand2005, Lele2008a}. In 1D, the model consists of the following macroscopic and microscopic force balances written in that order:
\begin{gather}
	\tau_{,y} = 0, \label{eq:macro}\\
    \tau = \tau^p - k^p_{,y},\label{eq:micro}
\end{gather}
Consistent with the principles of rational thermodynamics, we also have the following constitutive equations.
\begin{gather}
	\tau = \mu(u_{,y} - \gamma^p), \\
	\tau^p = S \left(\frac{d^p}{d_0}\right)^m \frac{\dot{\gamma^p}}{d^p}, \label{eq:micro_consti}\\
	k^p = S_0 L^2 \gamma^p_{,y} + S_0 l^2 \left(\frac{d^p}{d_0}\right)^m \frac{\dot{\gamma^p}_{,y}}{d^p}, \label{eq:micro_consti2}\\
	\dot{S} = H(S)d^p, \quad S(y, 0) = S_0 > 0, \\
	d^p = \sqrt{|\dot{\gamma^p}|^2 + l^2 |\dot{\gamma^p}_{,y}|^2},
\end{gather}
where $\tau$ is the shear stress, $\tau^p$ and $k^p$ are respectively microstress and gradient microstress, $\mu$ is elastic shear modulus, $u$ is the displacement, $\gamma^p$ is the plastic strain, $S_0$ is the coarse-grain yield strength, $S$ is an internal state variable, $d^p$ is an effective flow-rate, $d_0$ is a reference flow-rate, $m<1$ is a rate sensitivity parameter, $L$ and $l$ are energetic and dissipative length-scales, and $H$ is a hardening function.

The model is solved subject to the following boundary conditions
\begin{gather}
	u(0, t) = 0, \quad u(h, t) = u^\dagger(t) \quad \text{(prescribed)}, \\
	\dot{\gamma^p}(0, t) = \dot{\gamma^p}(h, t) = 0, \label{eq:plastic_bc}
\end{gather}
and the initial conditions
\begin{equation}
	u(y, 0) = 0, \quad \gamma^p(y, 0) = 0. \label{eq:init}
\end{equation}

In 2D, the model consists of the following macroscopic and microscopic force balances written in that order:
\begin{gather}
    \text{div}\bm{T} + \bm{b} = 0, \label{eq:macro_2d}\\
    \text{div}\bm{\xi} + \tau - \pi = 0, \label{eq:micro_2d}
\end{gather}
The constitutive equations are
\begin{align}
    \bm{T} &= 2\mu (\bm{E} - \bm{E^p}) + \lambda(\text{tr}\bm{E})\bm{1},  \\
    \bm{\dot{E^p}} &= \dot{\gamma^p} \bm{N^p}, \quad \bm{N^p} = \frac{\bm{T_0}}{|\bm{T_0}|}, \quad \dot{\gamma^p} \geq 0, \\
    \gamma^p &= \int_0^t \dot{\gamma^p}(\zeta)d\zeta, \\
    \phi^p &= \int_0^t \dot{\eta^p}(\zeta)d\zeta, \quad \text{where} \quad \dot{\eta^p} = |(\nabla\dot{\gamma^p}\times)\bm{N^p}|, \\
    d^p &= l_3|\nabla\dot{\gamma^p}|,\\
    \tau &= |\bm{T_0}|, \\
    \pi &= \left(S_0 \sqrt{(f(\gamma^p))^2 + l_2 \phi^p}\right) \left(\frac{\dot{\gamma^p}}{\nu_0}\right)^m, \quad \text{where} \quad f(\gamma^p) = 1 + \frac{H_0}{S_0}\gamma^p, \label{eq:micro_consti2d}\\
    \bm{\xi} &= S_0 l_1^2\nabla\gamma^p + S0\left(\frac{d^p}{d_0}\right)^q l_3^2 \frac{\nabla\dot{\gamma^p}}{d^p}, \label{eq:micro_consti2d_2}
\end{align}
where $\bm{T}$ is the Cauchy stress, $\bm{b}$ the body force, $\bm{\xi}$ and $\pi$ microstresses, $\bm{E^p}$ the plastic part of the total strain $\bm{E}$, $\dot{\gamma^p}$ the equivalent plastic strain rate, $\gamma^p$ the equivalent plastic strain, $\bm{N^p}$ the plastic flow direction, $\bm{T_0}$ the deviatoric stress, $\tau$ the deviatoric stress resolved in the flow direction, $f(\gamma^p)$ a linear isotropic hardening function, $l_1$, $l_2$, and $l_3$ the length-scales and $\mu$ and $\lambda$ the elastic Lam\'e moduli. $S_0 > 0$ is the initial flow strength, $d_0$ and $\nu_0$ are the reference flow rates, and $m$ and $q$ are rate sensitivity parameters.
The field equations are solved subject to the displacement boundary conditions
\begin{equation}
    \bm{u}(x_1, 0, t) = \bm{0}, \quad u_1(x_1, h, t) = u^*(t) \quad \text{(prescribed)}, \quad u_2(h, t) = 0;
\end{equation}
microscopically hard boundary conditions
\begin{equation}
    \dot{\gamma^p}(x_1, 0, t) = \dot{\gamma^p}(x_1, h, t) = 0;
\end{equation}
and initial conditions
\begin{equation}
    \bm{u}(x_1, x_2, 0) = \bm{0}, \quad \gamma^p(x_1, x_2, 0) = 0. \label{eq:init_2d}
\end{equation}

\begin{figure}
    \centering
    \includegraphics[width=0.5\linewidth]{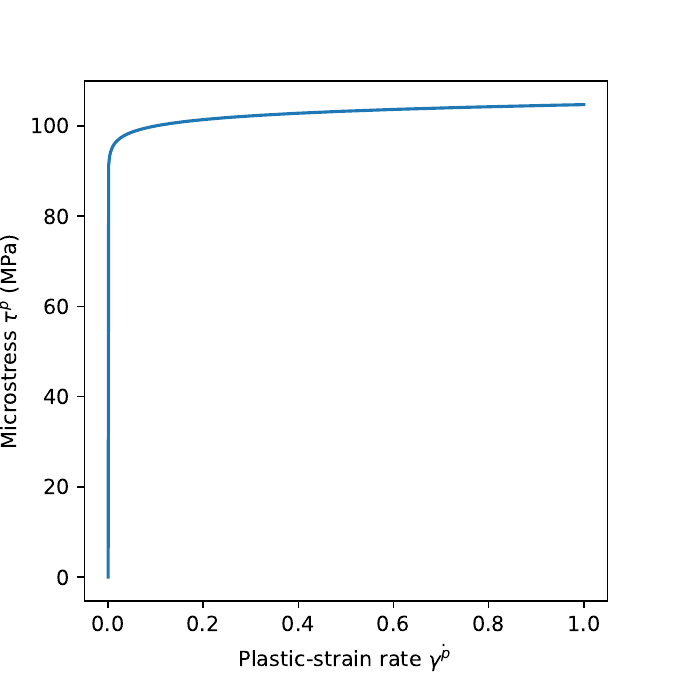}
    \caption{Evolution of microstress with plastic-strain rate, while $H=l=0$, $S_0=100$ MPa, $m=0.02$, and $d0=0.1$.}
    \label{fig:microstress}
\end{figure}

Note that these balances as a whole or even the microforce balance considered separately (equations~\ref{eq:micro} and \ref{eq:micro_2d}) are very stiff
thanks to the nature of constitutive equations of microstresses (equations~\ref{eq:micro_consti}, \ref{eq:micro_consti2}, \ref{eq:micro_consti2d}, and \ref{eq:micro_consti2d_2}). For example, $\tau^p$ in equation~\ref{eq:micro_consti} is not differentiable in $d^p$ at $\gamma^p=0$ as $m\ll 1$. Ill conditioning is also due to the presence of length-scales in the microforce balance. These length scales are small relative to the typical stiffness coefficients in the macroforce balance for metals.  As $\gamma^p$ and its derivatives are close to zero at the beginning of simulations, the derivatives of the microforce-balance residue would be very large, precipitating numerical difficulties. To illustrate this, we plot the microstress $\tau^p$ (equation~\ref{eq:micro_consti}) whilst keeping the hardening function $H$ and dissipative length-scale $l$ zero; see figure~\ref{fig:microstress}.
 In the figure one observes a vertical rise of the microstress followed by a very sharp transition to a gradual rise. This extremely stiff behavior of the microstress is also reflected in the stress-strain curves, frequently causing convergence failures in the finite element (FE) implementation. In fact, in most FE implementations, the transition regions are completely avoided by using large time-steps, and the complete stress-strain curve is obtained through back-tracing and extrapolation~\cite{Anand2005, Lele2008a, Lele2008}. Even with such bypasses, some special measures are needed to make the Newton-Raphson scheme converge. The authors of this article have not been successful in the FE implementation of the theory due to convergence failures, and an existing implementation reported in~\cite{Lele2008} fails to run. Accordingly, we will be comparing our PIDL simulation plots with digitized versions of the plots given in the literature~\cite{Anand2005}.
Yet another drawback of the currently adopted SGP model is that it develops plasticity even while unloading. To see this clearly, one switches off the length-scales ($L=l=0$) in the one-dimensional model. Then the microscopic force balance reduces to $\tau = \tau^p$. While unloading, $\tau$ is non-zero and so is $\tau^p$. But this leads to a non-trivial expression of the plastic-strain rate -- see equation~\ref{eq:micro_consti}. Indeed, one continues to detect this anomaly even when the length-scales are nonzero.

Presence of the length scales (which impart a multi-scale structure to the model) and the piecewise polynomial nature of the approximants are additional reasons behind the general inadequacy of a $C^0$ FE scheme. Unlike an FE scheme, PIDL uses a non-incremental approach with neural networks as global (not piecewise) approximants that are, with a proper choice of the activation function, $C^\infty$. PIDL thus ameliorates several shortcomings of FE methods and provides consequent advantages. We accordingly expect PIDL methods to provide a more streamlined and better performing numerical framework for such problems as the SGP. 

\section{PIDL implementation and results}
The SGP models described in the last section and summarised in equations~\ref{eq:macro} to \ref{eq:init} and equations~\ref{eq:macro_2d} to \ref{eq:init_2d} are solved using PIDL methods with different values of length-scales and hardening modulus. With its generality, PIDL affords choices at every step of the implementation. To approximate the field variables, there are different neural networks to choose from, e.g. deep feedforward networks, convolutional networks, and recurrent networks~\cite{Goodfellow2016}. In this work, we use deep feedforward networks (DFNs), which are generally adopted for regression and classification tasks. The inputs to the network are the spatial and time coordinates. We may use a single network with all the field variables of interest as its outputs, or we may use one separate network for each field variable~\cite{Haghighat2021}. The first approach is more efficient as the same learnt features are utilized to predict all the field variables. In the second approach, the same features might be learnt in  different networks. Even so, the second approach provides a more flexible architecture. Our experience is that both the approaches may yield satisfactory results. While the first requires a sufficiently large network, the second demands longer training. All the numerical results, presented later in the section, are obtained with a single network using tanh as the activation function.

Another choice is that of a global optimization scheme to find the network parameters. Because of various reasons, first-order gradient-based optimization algorithms (possibly with the provision for a perturbation to facilitate global search) are the norm in deep learning. Among these, we choose Adam~\cite{Kingma2014}, which is one of the most popular algorithms. Adam has many advantages, e.g., intuitive hyper-parameters and automatic step-size annealing. The learning rate provided to Adam acts as an upper bound to the step taken in the parameter space, and the step is also invariant to the scale of gradients, i.e. the absolute scale of the loss function does not matter. These characteristics of Adam provide helpful guidance during the loss function formulation and hyper-parameter selection. For example, the scaling of the total loss value is not required, and one may choose the learning rate depending on how much change is allowed in the network parameters in a single update step. Although Adam has an inbuilt step-size annealing, one may enforce a finer control over the training by externally providing an annealing scheme for the learning rate.  We have used a linearly decreasing learning rate, which starts at $0.01$ and ends at $0.00001$.

As an aid to the optimization process, the inputs and outputs to the network should be of order one. We scale our inputs and outputs accordingly. The spatial and time coordinates $x$ and $t$ are thus divided by their maximum values. The displacement is divided by the maximum applied displacement loading, and the plastic strain is divided by $S_0/\mu$. To further simplify the optimization, we impose the boundary and initial conditions strictly by encoding them into the network architecture~\cite{Lu2021}.
 
 Our loss function is the average of squared residues of the field equations at the collocation points. The collocation points, also called the training points, are sampled uniformly over the problem domain $\Omega \times [t_{min}, t_{max}]$, where $\Omega$ denotes the spatial domain of the problem. We also sample a set of validation points, different from the training points, for which we monitor the loss value. We stop the training if the training loss is sufficiently low, or if the validation loss starts rising, which is a sign of overfitting. After the training, we plot the distributions of field variables over a grid of points to see how well the model predicts over the problem domain. We also plot the stress-strain curves and the equivalent plastic strain profiles, and compare them with results digitized from the available literature~\cite{Anand2005}; the digitized plots are acronymed DAL in the figures. The implemented programs are available at \url{https://github.com/Debasish-Roy-IISc/PIDL-gradient-plasticity}.

\subsection{Results of forward simulations}
\begin{figure}
     \centering
     \begin{subfigure}[b]{0.32\textwidth}
         \centering
         \includegraphics[width=\textwidth]{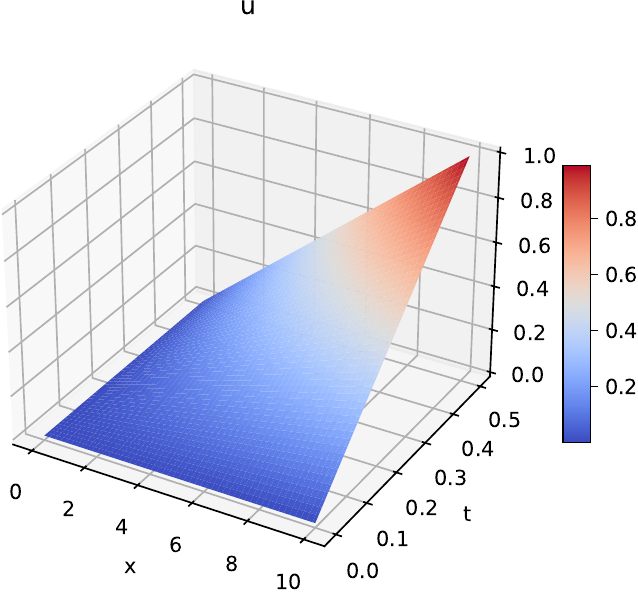}
         \caption{Predicted displacement}
     \end{subfigure}
     \hfill
     \begin{subfigure}[b]{0.32\textwidth}
         \centering
         \includegraphics[width=\textwidth]{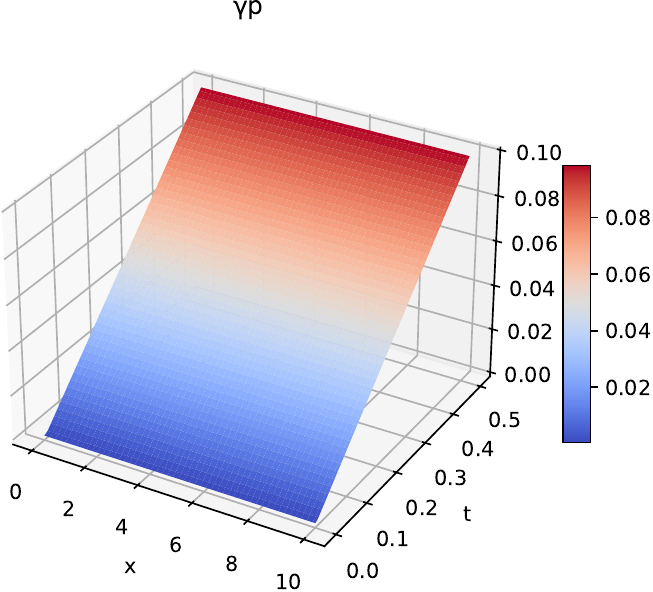}
         \caption{Predicted plastic strain}
     \end{subfigure}
     \hfill
     \begin{subfigure}[b]{0.32\textwidth}
         \centering
         \includegraphics[width=\textwidth]{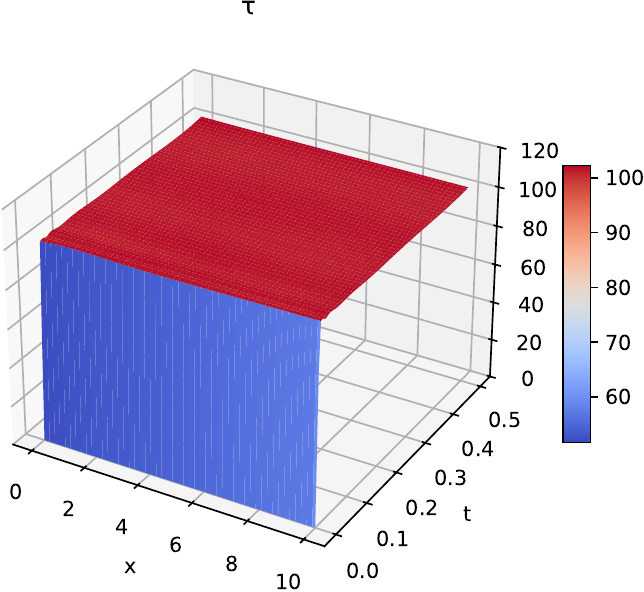}
         \caption{Predicted shear stress}
     \end{subfigure}
      \begin{subfigure}[b]{0.49\textwidth}
         \centering
         \includegraphics[width=\textwidth]{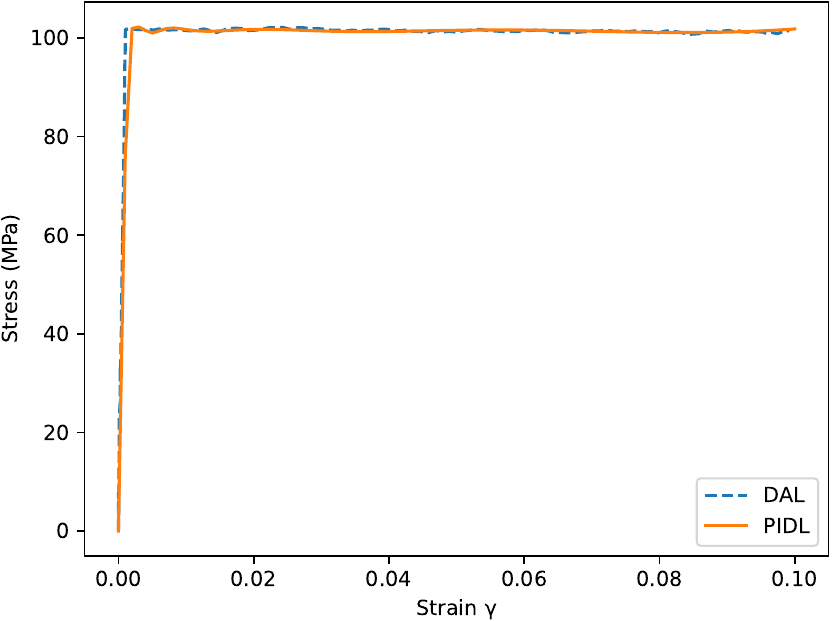}
         \caption{Stress-strain curves}
     \end{subfigure}
     \hfill
   \begin{subfigure}[b]{0.49\textwidth}
         \centering
         \includegraphics[width=\textwidth]{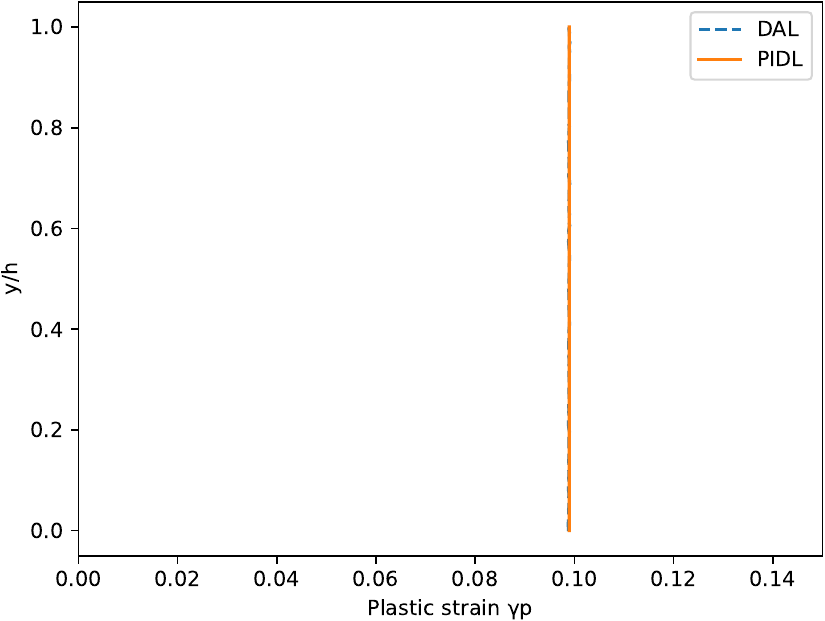}
         \caption{Plastic-strain profile across the strip}
     \end{subfigure}
        \caption{Predictions for no energetic and dissipative gradient hardening, and no internal variable hardening.}
        \label{fig:le0_ld0_H0}
\end{figure}

Numerical results for the 1D model shown below use the following material parameters:
\begin{equation}
    \mu = 100 \text{ GPa}, \quad S_0 = 100 \text{ MPa}, \quad d_0 = 0.1 \text{ s\textsuperscript{-1}}, \quad m = 0.02.
\end{equation}
Only linear internal-variable hardening is presently considered:
\begin{equation}
    H(S) = H \equiv \text{constant.}
\end{equation}
Results for the case when there are no energetic and dissipative gradient hardening ($L = l = 0$) and no internal variable hardening ($H=0$), are shown in figure~\ref{fig:le0_ld0_H0}. With this parameter setting, the model is close to (though not the same as) the rate-independent elastic perfectly-plastic theory. The SGP model is yield-free, and two field equations -- microscopic and macroscopic force balances -- must be solved simultaneously. Also, there is still hardening arising out of the power-law structure of the microforce constitutive equation~\ref{eq:micro_consti}. We have trained a network with three hidden layers, each having 32 neurons, for 25,000 epochs to obtain the results. The loss was minimized over 1,000 collocation points. The network outputs are the displacement $u$ and the plastic strain $\gamma^p$.

\begin{figure}
     \centering
     \begin{subfigure}[b]{0.32\textwidth}
         \centering
         \includegraphics[width=\textwidth]{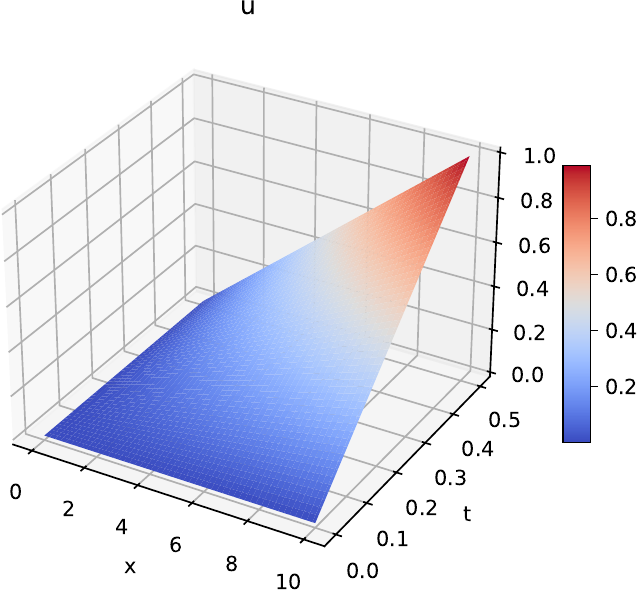}
         \caption{Predicted displacement}
     \end{subfigure}
     \hfill
     \begin{subfigure}[b]{0.32\textwidth}
         \centering
         \includegraphics[width=\textwidth]{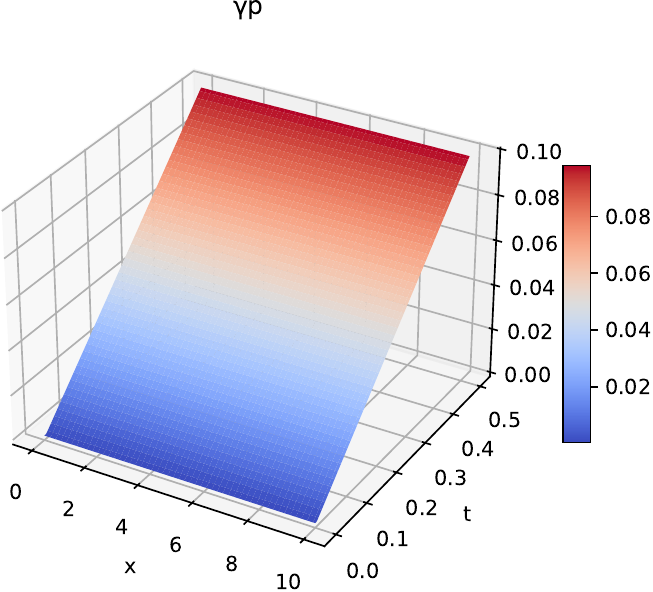}
         \caption{Predicted plastic strain}
     \end{subfigure}
     \hfill
     \begin{subfigure}[b]{0.32\textwidth}
         \centering
         \includegraphics[width=\textwidth]{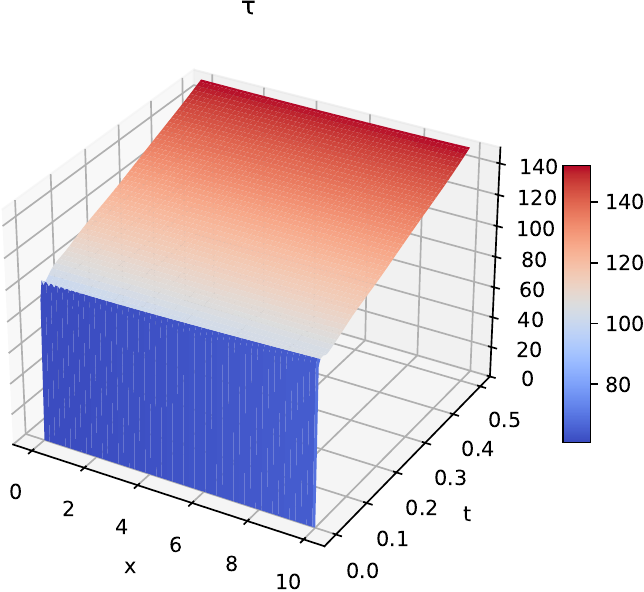}
         \caption{Predicted shear stress}
     \end{subfigure}
      \begin{subfigure}[b]{0.49\textwidth}
         \centering
         \includegraphics[width=\textwidth]{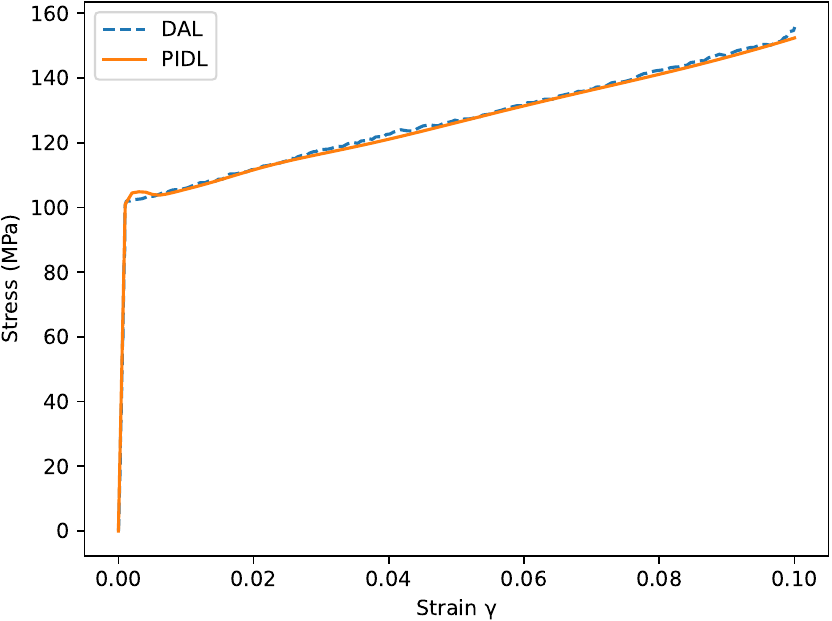}
         \caption{Stress-strain curves}
     \end{subfigure}
     \hfill
   \begin{subfigure}[b]{0.49\textwidth}
         \centering
         \includegraphics[width=\textwidth]{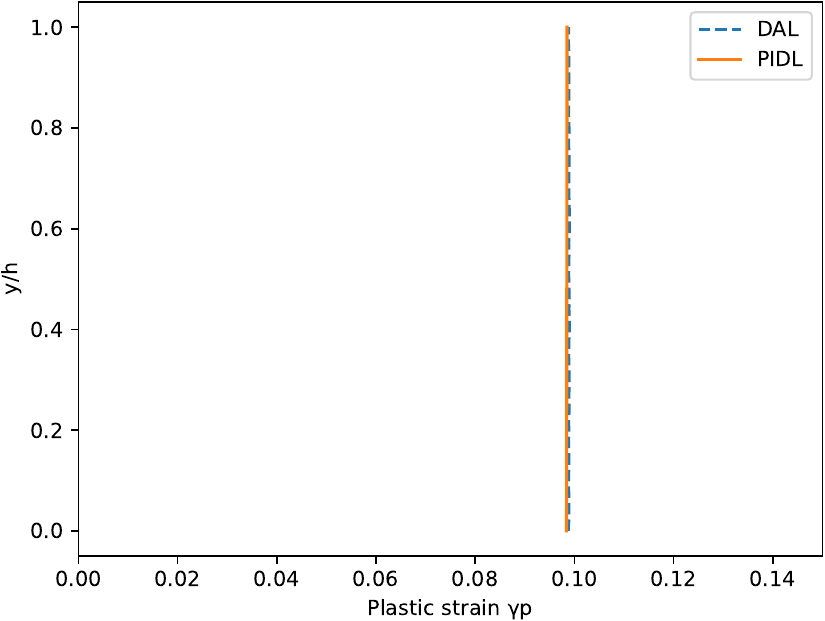}
         \caption{Plastic-strain profile across the strip}
     \end{subfigure}
        \caption{Predictions for internal variable hardening ($H = 500$ MPa)}
        \label{fig:le0_ld0_H500}
\end{figure}

\begin{figure}
     \centering
     \begin{subfigure}[b]{0.32\textwidth}
         \centering
         \includegraphics[width=\textwidth]{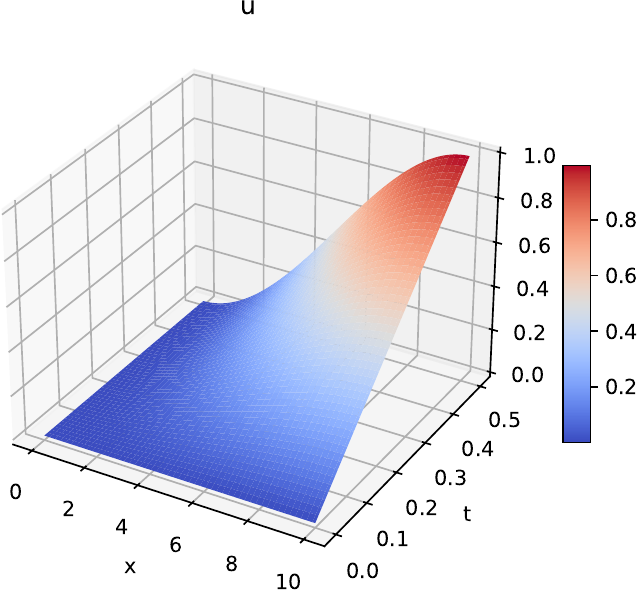}
         \caption{Predicted displacement}
     \end{subfigure}
     \hfill
     \begin{subfigure}[b]{0.32\textwidth}
         \centering
         \includegraphics[width=\textwidth]{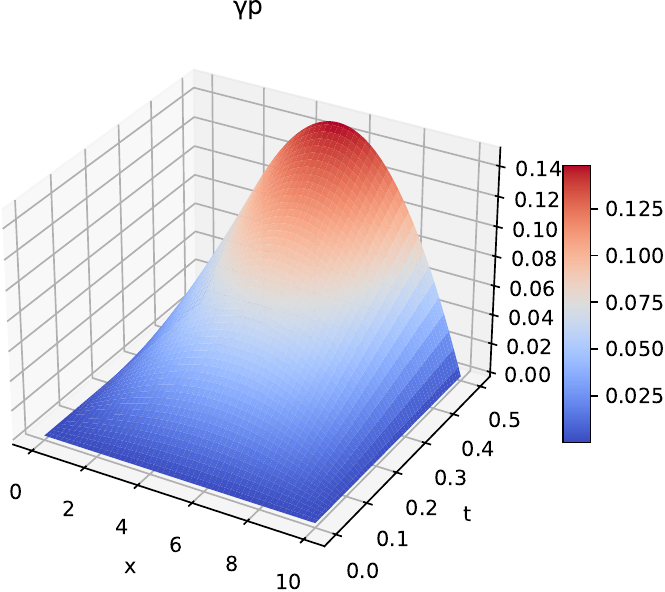}
         \caption{Predicted plastic strain}
     \end{subfigure}
     \hfill
     \begin{subfigure}[b]{0.32\textwidth}
         \centering
         \includegraphics[width=\textwidth]{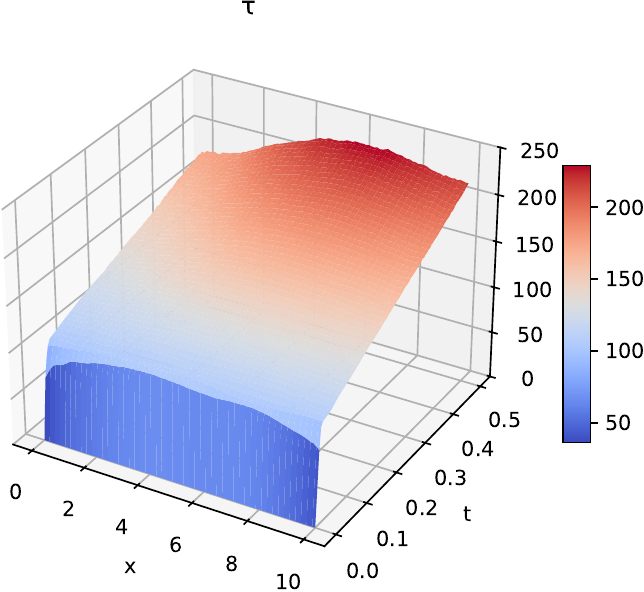}
         \caption{Predicted shear stress}
     \end{subfigure}
      \begin{subfigure}[b]{0.49\textwidth}
         \centering
         \includegraphics[width=\textwidth]{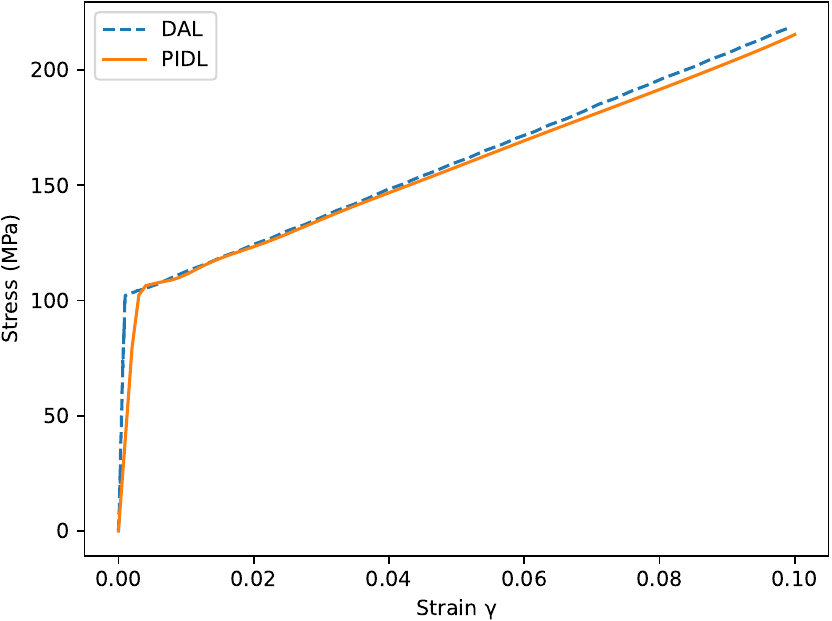}
         \caption{Stress-strain curves}
     \end{subfigure}
     \hfill
   \begin{subfigure}[b]{0.49\textwidth}
         \centering
         \includegraphics[width=\textwidth]{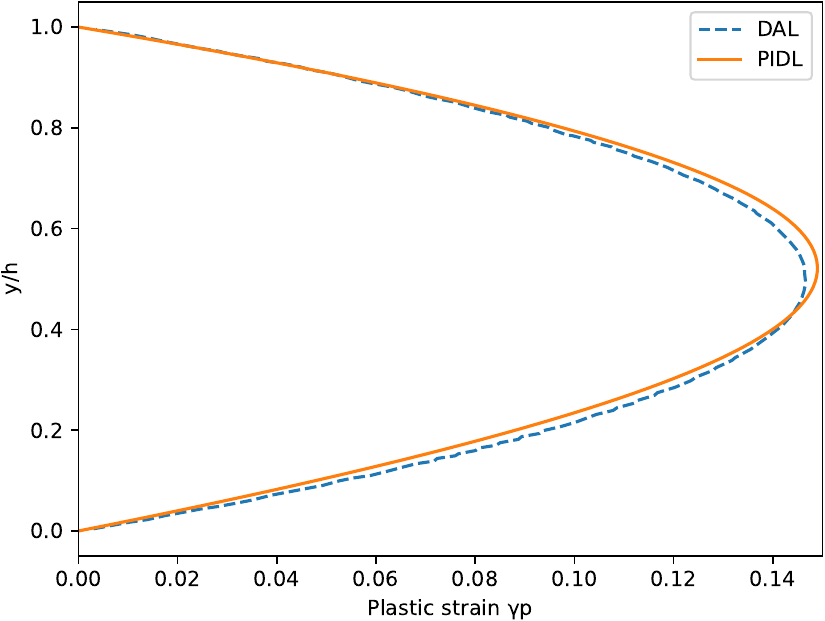}
         \caption{Plastic-strain profile across the strip}
     \end{subfigure}
        \caption{Predictions for nonzero energetic gradient hardening ($L = 10$)}
        \label{fig:le10_ld0_H0}
\end{figure}

We also consider the case when there is internal variable hardening ($H=500$ MPa) but no energetic or dissipative gradient hardening ($L = l = 0$); the results are shown in figure~\ref{fig:le0_ld0_H500}. Again, we have used a single network with three hidden layers, each with 32 neurons. The network outputs, as before, are $u$ and $\gamma^p$, and it is trained for just 25,000 epochs with 1,000 collocation points.
Results for the case when there is energetic-gradient hardening ($L = 10$), but no other hardening ($l = H = 0$) are provided in figure~\ref{fig:le10_ld0_H0}. For this case, we have utilized a network with three hidden layers, each having 64 neurons, to make it more expressive. We have also trained it for 50,000 epochs and used 10,000 collocation points (as a bigger network requires a longer training).

\begin{figure}
     \centering
     \begin{subfigure}[b]{0.32\textwidth}
         \centering
         \includegraphics[width=\textwidth]{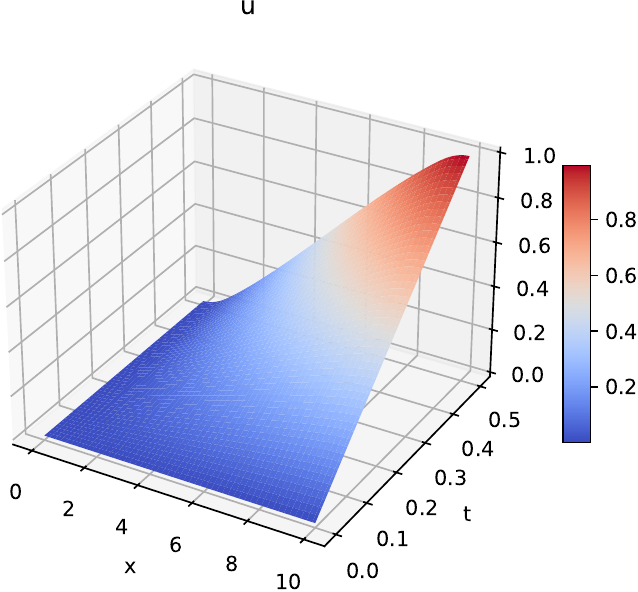}
         \caption{Predicted displacement}
     \end{subfigure}
     \hfill
     \begin{subfigure}[b]{0.32\textwidth}
         \centering
         \includegraphics[width=\textwidth]{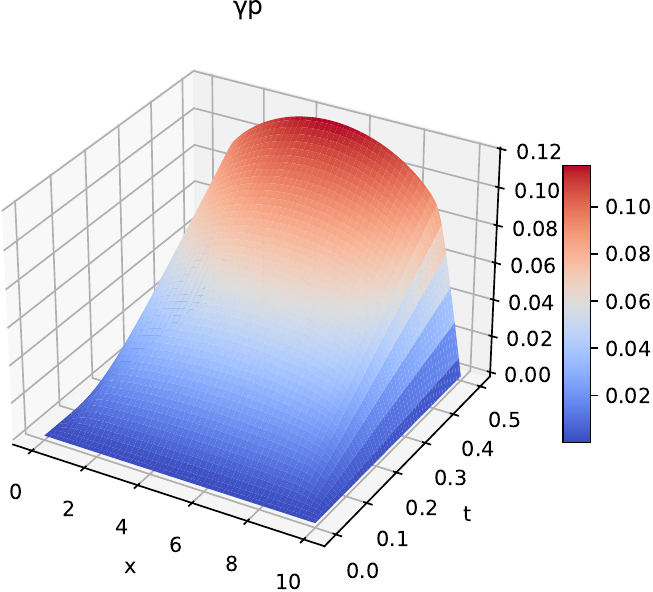}
         \caption{Predicted plastic strain}
     \end{subfigure}
     \hfill
     \begin{subfigure}[b]{0.32\textwidth}
         \centering
         \includegraphics[width=\textwidth]{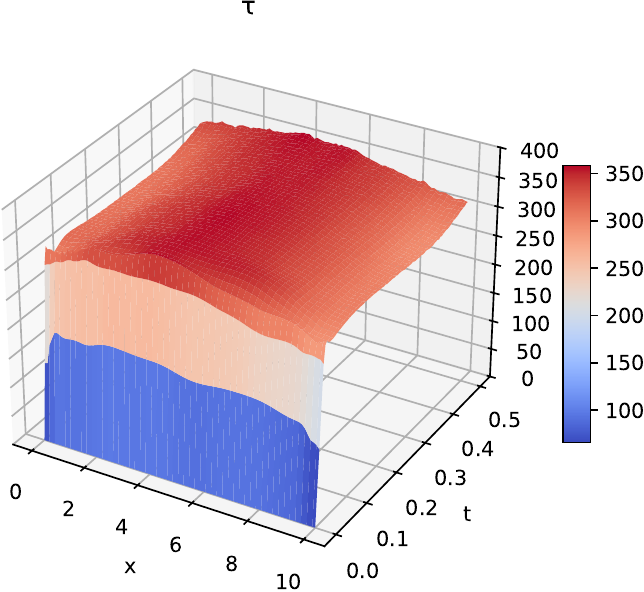}
         \caption{Predicted shear stress}
     \end{subfigure}
      \begin{subfigure}[b]{0.49\textwidth}
         \centering
         \includegraphics[width=\textwidth]{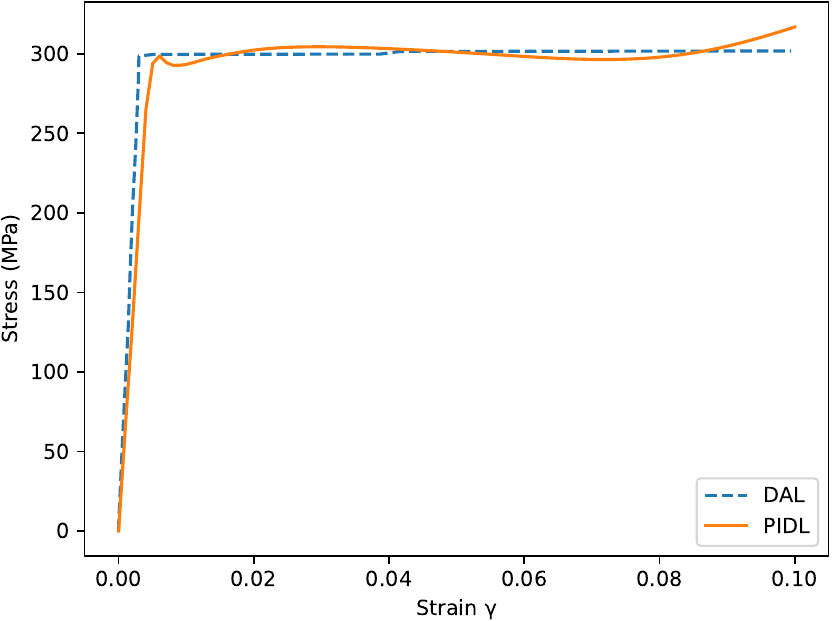}
         \caption{Stress-strain curves}
     \end{subfigure}
     \hfill
   \begin{subfigure}[b]{0.49\textwidth}
         \centering
         \includegraphics[width=\textwidth]{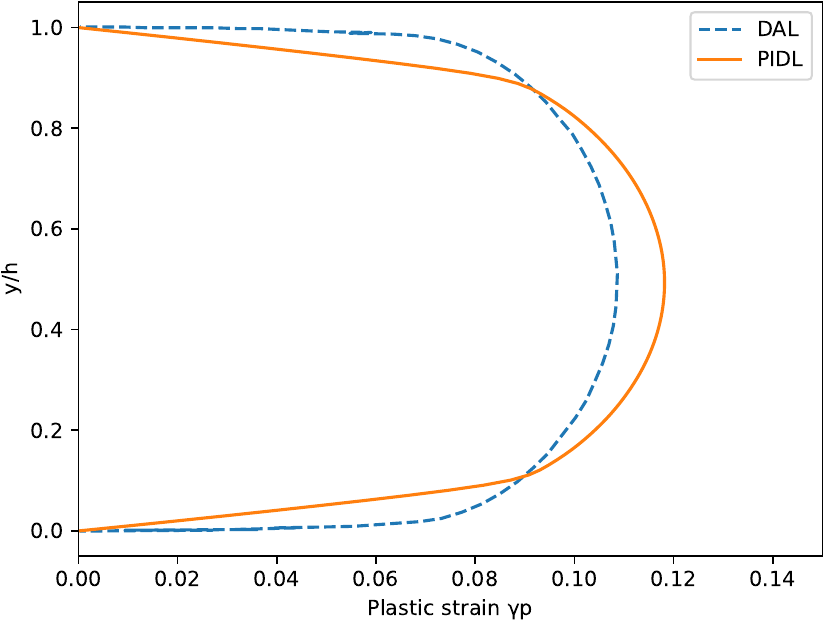}
         \caption{Plastic-strain profile across the strip}
     \end{subfigure}
        \caption{Predictions for nonzero dissipative-gradient hardening ($l = 10$)}
        \label{fig:le0_ld10_H0}
\end{figure}

 For a non-zero dissipative length-scale, the microscopic force balance is a fourth order PDE. To reduce the order of differentiation required, we use a mixed approach for this particular case. Specifically, the microgradient stress $k^p$ is also made a network output along with $u$ and $\gamma^p$. Again we use a network with three hidden layers, 64 neurons each, and train it over 10,000 collocation points for 100,000 epochs. The results are furnished in  figure~\ref{fig:le0_ld10_H0}.

\begin{figure}
     \centering
     \begin{subfigure}[b]{0.32\textwidth}
         \centering
         \includegraphics[width=\textwidth]{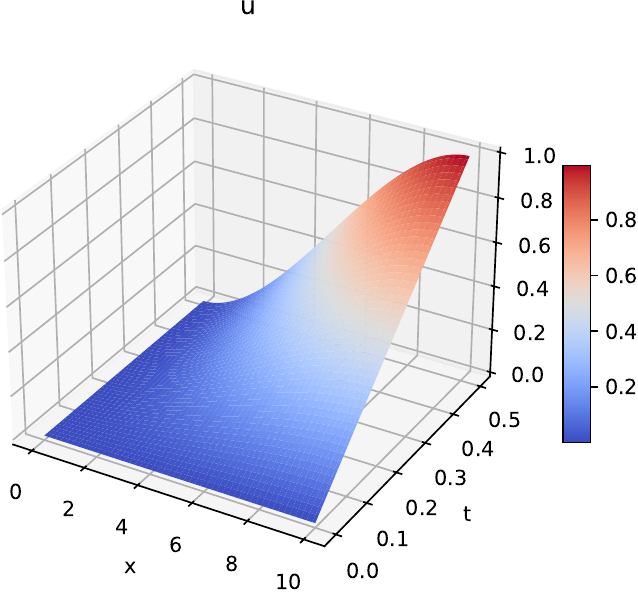}
         \caption{Predicted displacement}
     \end{subfigure}
     \hfill
     \begin{subfigure}[b]{0.32\textwidth}
         \centering
         \includegraphics[width=\textwidth]{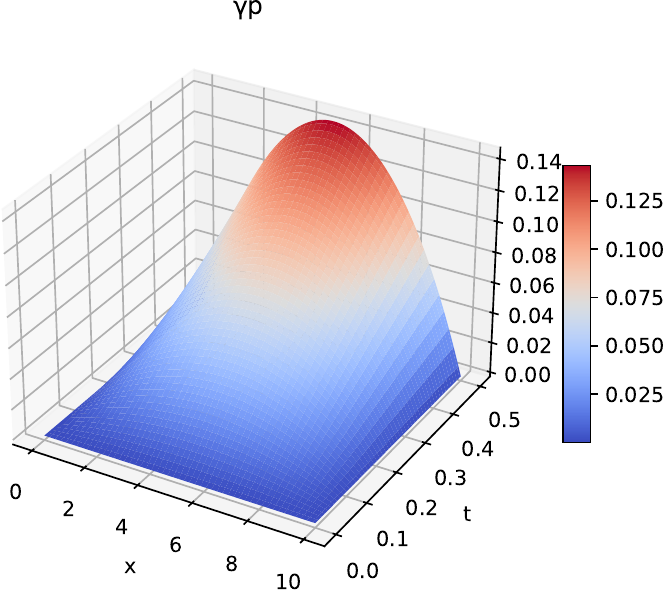}
         \caption{Predicted plastic strain}
     \end{subfigure}
     \hfill
     \begin{subfigure}[b]{0.32\textwidth}
         \centering
         \includegraphics[width=\textwidth]{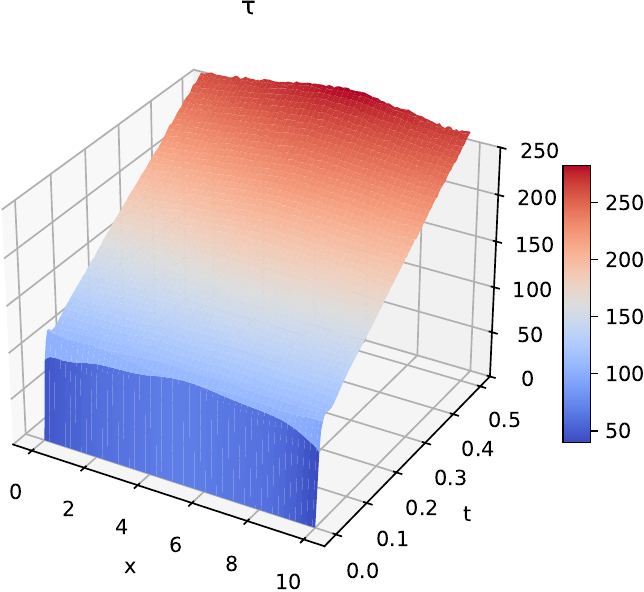}
         \caption{Predicted shear stress}
     \end{subfigure}
      \begin{subfigure}[b]{0.49\textwidth}
         \centering
         \includegraphics[width=\textwidth]{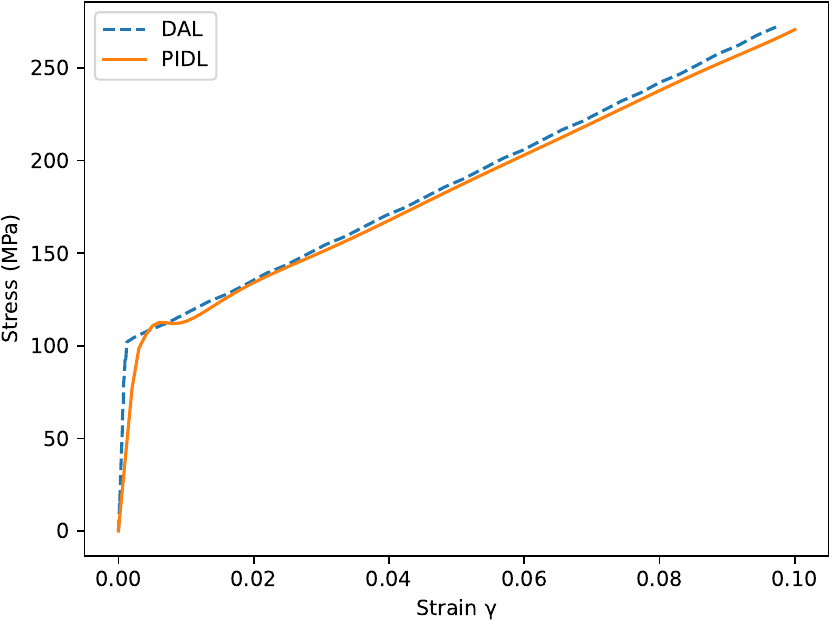}
         \caption{Stress-strain curves}
     \end{subfigure}
     \hfill
   \begin{subfigure}[b]{0.49\textwidth}
         \centering
         \includegraphics[width=\textwidth]{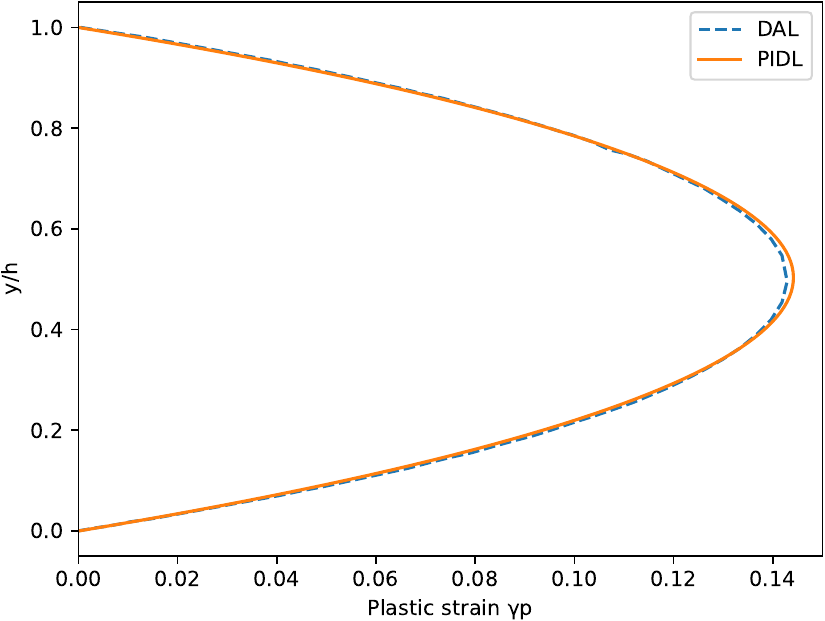}
         \caption{Plastic-strain profile across the strip}
     \end{subfigure}
        \caption{Predictions for combined energetic-gradient and internal variable hardening ($L = 10$ and $H = 500$ MPa)}
        \label{fig:le10_ld0_H500}
\end{figure}

Results for the case when there is combined energetic-gradient ($L=10$) and internal variable hardening ($H=500$ MPa) are shown in figure~\ref{fig:le10_ld0_H500}. The network used has three hidden layers, each with 64 neurons. The network is trained for 100,000 epochs with 10,000 collocation points.

FE implementation of the two-dimensional model suffer from issues similar to the ones already pointed out~\cite{Lele2008}. For one, the coupled field equations are stiff, leading to ill-conditioning and numerical corruption. There are a large number of terms in the Jacobian of residues, and a few terms that couple displacement and plastic strain are very expensive to evaluate. The model is typically made less stiff and computationally tractable by simply neglecting these terms~\cite{Lele2008}. These approximations nevertheless might be the cause of certain discrepancies in the FE-based simulation results. For instance, even though one attempts to solve the same physical system using 1D and 2D SGP models, the so-called elastic arm of the stress-strain curve (indeed, there is no purely elastic regime in a rate-dependent model) comes out far steeper with the 1D model. Moreover, in the 1D model, the sharp transition in the stress-strain curve occurs when the stress value equals the initial flow resistance ($S_0$). In contrast, this transition happens in the 2D model  at a stress value ($\approx 100$ MPa) below  $S_0=141.4$ MPa.

\begin{figure}
    \centering
    \begin{subfigure}[b]{0.49\textwidth}
        \centering
        \includegraphics[width=\linewidth]{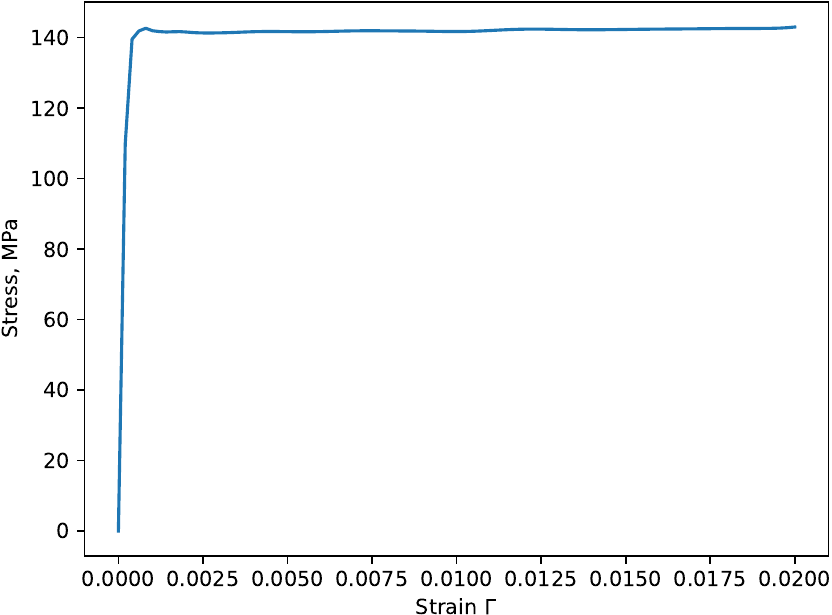}
        \caption{Stress-strain curve}
    \end{subfigure}
    \hfill
    \begin{subfigure}[b]{0.49\textwidth}
        \centering
        \includegraphics[width=\linewidth]{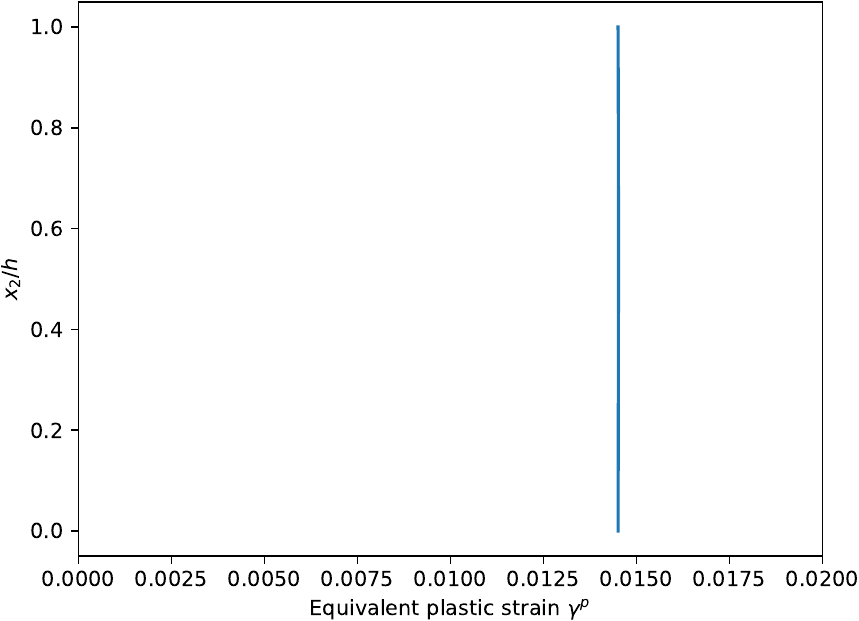}
        \caption{Plastic strain profile across the strip}
    \end{subfigure}
    \caption{PIDL simulation results for two-dimensional model with no hardening}
    \label{fig:2d_nohard}
\end{figure}

PIDL directly works on the strong form, neglecting no coupling terms. None of the discrepancies mentioned above appear in the PIDL simulation results. To illustrate this, we have used a PIDL implementation of the 2D model with no hardening, leading to results shown in figure~\ref{fig:2d_nohard}. For this simulation, we have used a DFN with three hidden layers, each having 64 neurons, and trained it for 30,000 epochs over 10,000 collocation points. The mechanical parameters used are as follows:
\begin{equation}
    E = 210 \text{ GPa}, \quad \nu = 0.3, \quad S_0 = 141.4 \text{ MPa}, \quad d_0 = 0.02828 \text{ s\textsuperscript{-1}}, \quad m = 0.05,
\end{equation}
where $E$ is the Young's modulus and $\nu$ the Poisson's ratio.

\subsection{Results of parameter interpolations}
Unlike the FE, PIDL models can perform acceptable forward simulations for parameter values on which they were never trained. To test this in the context of the SGP model, we treat the shear modulus ($\mu$) and the initial flow resistance ($S_0$) as inputs to the neural network (in addition to the space and time coordinates). For a ready demonstration, we train the model with no length-scales and no internal variable hardening for a range of values of $\mu$ and $S_0$. Once training is accomplished, we predict the stress-strain curves for other values of $\mu$ and $S_0$ not used in the training phase. Whenever possible, we compare the results with those available in the literature. 

\begin{figure}
    \centering
    \begin{subfigure}[b]{0.49\textwidth}
        \centering
        \includegraphics[width=\linewidth]{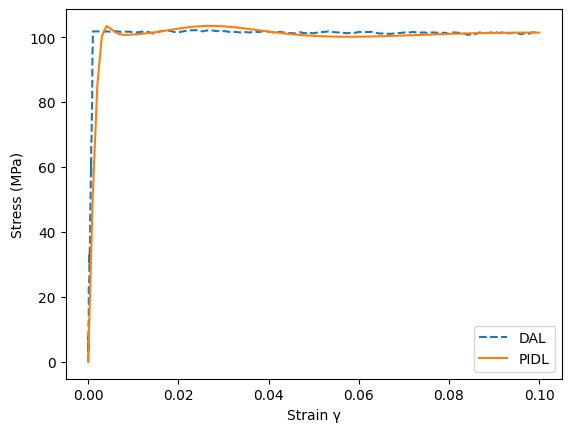}
        \caption{}
        \label{fig:mu_pred_a}
    \end{subfigure}
    \hfill
    \begin{subfigure}[b]{0.49\textwidth}
        \centering
        \includegraphics[width=\linewidth]{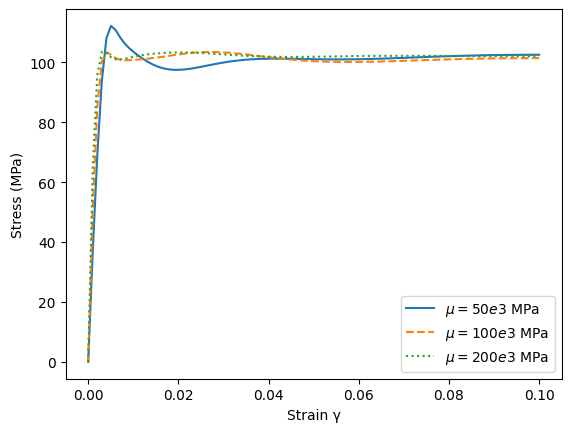}
        \caption{}
        \label{fig:mu_pred_b}
    \end{subfigure}
    \caption{Predictions of stress-strain curves for values of shear modulus $\mu$ on which the model was never trained; (a) a comparison with the result available in literature and (b) predictions for vastly different values of $\mu$.}
    \label{fig:mu_pred}
\end{figure}

\begin{figure}
    \centering
    \begin{subfigure}[b]{0.49\textwidth}
        \centering
        \includegraphics[width=\linewidth]{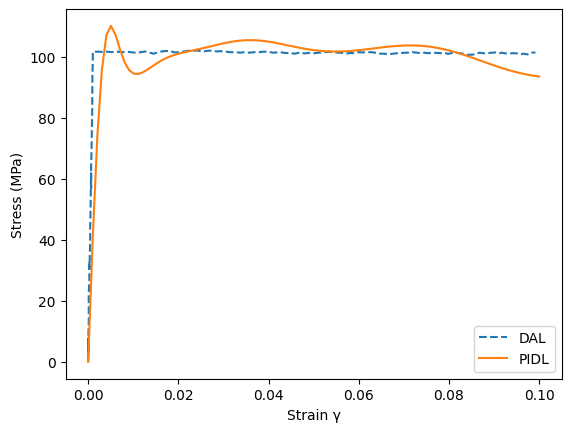}
        \caption{}
        \label{fig:S0_pred_a}
    \end{subfigure}
    \hfill
    \begin{subfigure}[b]{0.49\textwidth}
        \centering
        \includegraphics[width=\linewidth]{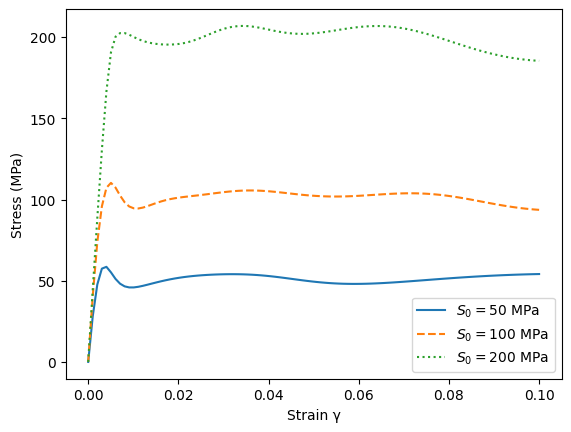}
        \caption{}
        \label{fig:S0_pred_b}
    \end{subfigure}
    \caption{Predictions of stress-strain curves for values of initial flow-resistance $S_0$ against which the model was never trained.}
    \label{fig:S0_pred}
\end{figure}

First, we train the model for a fixed $S_0$ (100 MPa) and a range of values of $\mu$ (30 equally spaced ones between 10 GPa and 1000 GPa). We predict the stress-strain curve for $\mu = 100$ GPa  and compare with results available in the literature (figure~\ref{fig:mu_pred_a}). The two curves match reasonably well. Note that the nearest training data point for $\mu = 100$ GPa is more than 12 GPa away. We also predict the stress-strain curves for vastly different values of $\mu$ (see figure~\ref{fig:mu_pred_b}). We expect $\mu$ to influence only the pre-transition portion of the stress-strain curve, the rise being steeper as the value of $\mu$ increases.  This is exactly what we see in figure~\ref{fig:mu_pred_b}.
Second, we train the model for a fixed value of $\mu$ (100 GPa) and a range of values of $S_0$ (30 equally spaced values between 10 MPa to 1000 MPa). Again, we predict the stress-strain curve for $S_0 = 100$ MPa and compare it with the one in the literature; see figure~\ref{fig:S0_pred_a}. We expect the point of transition in the stress-strain curve to shift upwards along the stress axis as we increase $S_0$. Our predictions are consistent with this expectation; see figure~\ref{fig:S0_pred_b}. Although, our predictions show some oscillations, we believe these might be removed through suitable modifications in the training.

We have experienced difficulties in training a network with both $\mu$ and $S_0$ as inputs. This may be due to the globally supported nature of the activation functions in the DFN. Indeed, the effects of $\mu$ and $S_0$ are very distinct and both have a spatiotemporally local character. Globally supported activation function may also be the cause of other issues we have faced with PIDL simulations, e.g. oscillations in the solution. A few of these shortcomings of a DFN-based PIDL are already recorded in the literature~\cite{wang2021understanding,wang2022respecting}; however the remedies recommended have not worked for us. If globally supported activation is one of the reasons that these difficulties arise, a changeover to neural networks with locally supported activation functions~\cite{Liu2024} might offer a consistent remedy.

\section{Conclusions}
Given the attractive functional approximation properties of deep neural networks~\cite{devore2021neural} and an inherent flexibility in the choice of network architectures, a natural question is whether such approximants could be exploited to design numerical schemes for stiff, multiscale and nonlinear problems in the mechanics of solids -- especially in cases where popular methods such as the FEM struggle. One such problem, considered in this article, is a strain gradient continuum model for metal plasticity wherein an accounting of the size effect (typically observed during plastic flows at micron scales) brings in length scales, thereby imparting to the model a multiscale structure. We have specifically used a physics-informed feedforward network where the training is based on a loss function that incorporates the underlying balance laws, constitutive relations and the boundary/initial conditions. The reported numerical results clearly attest to the feasibility of a neural approximation approach for the gradient plasticity model. With it comes the added advantages, e.g. the ability to extrapolate solutions to other parameters, load cases, initial conditions etc.

On the other hand, one may also observe limitations of the physics-informed network in its current form. For instance, it is not immediately clear how one could simulate cyclic or random loading with physics informed DFNs in the presence of inelasticity in the response. Also, the encoding of boundary and initial conditions, which simplify the optimization problem, may not always be possible. Another limitation is that the geometry of the problem is hard-coded into the network parameters during the learning process -- the network will not generalize over the geometry of the problem. By using graph neural networks (GNNs), we may be able to address these issues. A GNN takes the geometry and state of the problem as inputs and predicts the future geometry and state. Such a scenario is thus well suited to the simulations of deforming solids, where the base space (solid geometry) continues to evolve. Clearly, it also has the potential to extrapolate over initial domain geometries. Yet another issue with the DFN is its use of globally supported activation functions, such as tanh. As optimization happens for the entire domain simultaneously, any local singularity or ill-conditioning may have global effects -- possibly rendering the final solution sensitive to the parameters. It appears that the recently proposed Kolmogorov-Arnold Networks (KANs)~\cite{Liu2024}, which use B-splines as the activation functions, could aid in addressing this issue. We plan to take up GNNs and KANs in our future work in the context of nonlinear and inelasticity problems in the mechanics of solids, especially as a superior alternative to traditional methods, such as the FEM and peridynamics, SPH or other mesh-free methods.

\bibliography{refs.bib} \bibliographystyle{unsrt}

\end{document}